\documentclass{article}

\usepackage{graphics}
\usepackage{natbib}
\usepackage[colorlinks=true]{hyperref}
\bibpunct{(}{)}{;}{x}{,}{,}

\newcommand{\paren}[1]{\left ( {#1} \right )}
\newcommand{\ip}[2]{\left < {#1},{#2} \right >}

\newcommand{\M}{\mathcal{M}}
\newcommand{\X}{\mathcal{X}}

\newcommand{\ko}{\kappa_0}
\newcommand{\kt}{\kappa_2}

\newcommand{\matr}[1]{ \mathbf{#1} }

\begin{document}

\title{The Volume-of-Tubes formula: Computational Methods and Statistical
  Applications}
\author{
Catherine Loader \\
Department of Statistics \\
Case Western Reserve University \\
Cleveland, OH 44106
}

\maketitle

\abstract{

The volume-of-tube formula was first introduced by
\cite{hotelling}, to solve significance of terms in nonlinear regression
models. Since this pioneering paper, there has been significant work
on extending the tube formula to more general settings, including
multidimensional problems, and many new applications in statistical
inference, including confidence bands in regression and smoothing models;
applications to functional data analysis; testing in mixture models;
and spatial scan analysis.

Implementation of the tube formula requires numerical evaluation of certain
problem-specific geometric constants that appear in Hotelling's formula
and its extensions. The purpose of this note is to describe a software
library, \texttt{libtube}, that performs the calculations. A variety
of illustrative examples are given.

Source code for the \texttt{libtube} library and
examples can be downloaded from \url{http://www.herine.net/stat/libtube/}.
}

\section{Introduction}

The volume-of-tube problem can be stated rather simply. Given
a curve (or manifold) \(\M\) lying in \(n\)-dimensional Euclidean space,
what is the volume of the set of all points lying within a radius \(r\)
of the curve? In statistical applications, the spherical version of
this problem often arises; the manifold lies on the surface of the
unit sphere in \(n\) dimensions, and one wishes to compute the
\((n-1)\)-dimensional volume (or surface area) of the set of points
lying within a distance \(r\) of the manifold.

The volume-of-tubes formula was formulated and solved by
\citet{hotelling}, motivated
by application to significance testing in nonparametric regression.
A companion paper, \citet{weyl}, extended the results to higher dimensional
manifolds; that is, when \(\M\) is a surface, or more generally when
\(\M\) is a manifold of dimension \(d \le n\).

The main purpose of this article is to describe a set of routines
written by the author to implement the volume-of-tube formula in
statistical problems. In section \ref{sec:tube} the tube formula
(with boundary corrections) is described. The \texttt{libtube} software
is described in Section \ref{sec:libtube}. Applications to
non-linear regression, simultaneous confidence bands and
mixture modeling are described in Sections \ref{sec:nlreg},
\ref{sec:scb} and \ref{sec:mixture} respectively.

\section{The Volume-of-Tubes Formula}
\label{sec:tube}

\begin{figure}
\centerline{\scalebox{0.8}{\includegraphics{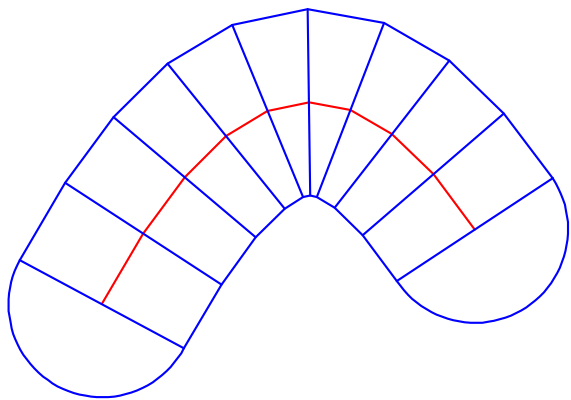}}}
\caption{The manifold is represented by the red curve. The tubular
  neighborhood is approximated by trapezoids, plus the two end-point caps.}
\label{fig:plotta}
\end{figure}

The volume-of-tubes formula was first derived by \cite{hotelling}.
The result can be illustrated on the plane by reference to Figure
\ref{fig:plotta}.
The manifold is represented by the red curve. The tubular
neighborhood of a given radius \(r\) is approximated by trapezoids,
plus the two end-point caps.
Adding up the area of the trapezoids and letting the partition become
increasingly fine shows that the area (or two-dimensional volume) of the
tube is
\[
  \hbox{Length of Manifold} \times 2r + \pi r^2.
\]

The \(2r\) represents the cross-sectional area of the manifold, while
\(\pi r^2\) represents the area of the end-point caps. The result
extends to manifolds embedded in \(n\)-dimensional space;
\[
  \hbox{Volume} = \kappa_0 V_{n-1} r^{n-1}  + \frac{l_0}{2} V_n r^n.
\]
Here \(\kappa_0\) is the length of the manifold and \(l_0\) is the number
of end-points (often, \(l_0=2\)). The functions \(\psi_0(r)\) are the
cross-sectional area and volume of the end-point caps respectively,
and \(V_k = \pi^{k/2}/ \Gamma( 1+k/2)\) is the volume of the \(k\)-dimensional
unit sphere.

When the manifold lies on the unit sphere, the result is similar, but the
cross-sectional area is replaced by a certain partial beta function. The result
is
\[
  \hbox{Volume} = \frac{\kappa_0 A_n}{2\pi} P(B_{1,(n-2)/2} \ge w^2)
    + \frac{l_0 A_n}{4} P(B_{1/2,(n-1)/2} \ge w^2),
\]
where \(B_{a,b}\) denotes a random variable following a beta distribution
with parameters \(a\) and \(b\); \(A_n = 2\pi^{n/2}/\Gamma(n/2)\)
is the surface area of the unit sphere in \(R^n\), and \(w = 1-r^2/2\).

\emph{Multidimensional Manifolds}

\begin{figure}
\centerline{\scalebox{0.8}{\includegraphics{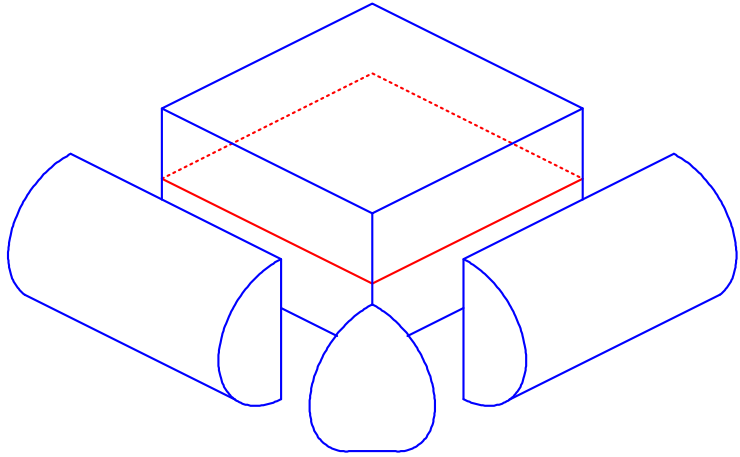}}}
\caption{Tube around a two dimensional manifold. The manifold is
  shown in red, and the tube is divided into a main part,
  half-cylinders around the edges, and corner wedges.}
\label{fig:bivman}
\end{figure}

Figure \ref{fig:bivman} shows a tube around a two-dimensional manifold.
To compute the volume of the tubular neighborhood, one divides the tube
into different pieces: a main piece, the half-cylinders around each
edge, and wedges at each corner of the manifold. For higher dimensional
manifolds, the ideas are similar, but there are more pieces to take care of.

A version of the tube formula, without boundary corrections, was first
derived by \cite{weyl}. \cite{naiman90} provided boundary corrections.
The result is a series with \(d+1\) terms. The first four terms
(for manifolds on the unit sphere) are
\begin{eqnarray*}
  \hbox{Volume} &=&
    \frac{ \ko A_n}{ A_{d+1} } P( B_{(d+1)/2,(n-d-1)/2} > w^2 ) \\
  && \; + \frac{ l_0 A_n}{ 2 A_d } P( B_{d/2,(n-d)/2} > w^2 ) \\
  && \; + \frac{\kt + l_1 + m_0}{2\pi} \frac{A_n}{A_{d-1}}
      P( B_{(d-1)/2,(n-d+1)/2} > w^2 ) \\
  && \; + \frac{l_2 + m_1 + n_0}{4\pi} \frac{A_n}{A_{d-2}}
      P( B_{(d-2)/2,(n-d+2)/2} > w^2 ).
\end{eqnarray*}
The constants \(l_0\), \(l_1\) and \(l_2\) arise from the corresponding series
for the half-tubes around the boundaries of the manifold.
 \(l_0\) is the \((d-1)\)-dimensional volume of the
boundaries (or the total length of the four edges in Figure \ref{fig:bivman}).
\(l_1\) and \(l_2\) are higher order terms representing boundary curvature.

\(m_0\) and \(m_1\) arise from the `corner wedges' where two
boundary faces meet. In figure \ref{fig:bivman}, \(m_0\) is the
sum of the four wedge angles at each corner of the manifold.

\(n_0\) arises for manifolds with \(d \ge 3\), from the corners where
three (or more) boundary faces meet.

\subsection{Random Processes}

In statistical applications, the fundamental use of the tube formula is
to find (or at least approximate) the distribution of the maximum of certain
random processes. As a simple example, consider the process
\[
  Z(\lambda) = \ip{T(\lambda)}{U}
\]
where \(T(\lambda)\) is an \(R^n\)-valued vector function, and \(U\) is
uniformly distributed over the unit sphere. Suppose that one is interested in
finding
\[
  P(\sup_{\lambda} Z(\lambda) \ge w) \label{eq:uniprob}
\]
for some \(w\).

The inner product exceeds \(w\) if, and only if, \(U\) is sufficiently close
to \(T(\lambda)\). Specifically,
\[
  \|T(\lambda) - U\|^2 = \|T(\lambda)\|^2 - 2 \ip{T(\lambda)}{U} + \|U\|^2
    = 2( 1- \ip{T(\lambda)}{U} ).
\]
Hence, \( \|T(\lambda) - U\| \le r \) if, and only if,
\( \ip{T(\lambda)}{U} \ge w \), where \( r^2 = 2(1-w)\). The probability
(\ref{eq:uniprob}) is therefore simply
\[
  \frac{ \hbox{ Area of tube of radius \(r\) around  \(\{T(\lambda)\}\)} }{
    \hbox{Surface area of unit sphere in \(R^n\)} }.
\]

In many statistical applications, one is interested in the distribution
of the maximum of a Gaussian process,
\[
  Z(\lambda) = \ip{ T(\lambda)}{\epsilon}
\]
where \(\epsilon\) follows the standard multivariate normal distribution.
To reduce this to the uniform process, one needs to condition on the
length of the \(\epsilon\) vector, and integrate over the conditional
distribution; see \cite{sunloascb}. The final result, up to fourth order, is
\begin{eqnarray*}
  P( \sup Z(\lambda) \ge c) &\approx&
    \frac{ \ko }{ A_{d+1} } P( \chi^2_{(d+1)/2} \ge c^2 ) \\
  && \; + \frac{ l_0 }{ 2 A_d } P( \chi^2_{d/2} > c^2 ) \\
  && \; + \frac{\kt + l_1 + m_0}{2\pi A_{d-1}} P( \chi^2_{(d-1)/2} > c^2 ) \\
  && \; + \frac{l_2 + m_1 + n_0}{4\pi A_{d-2}} P( \chi^2_{(d-2)/2} > c^2 ),
\end{eqnarray*}
where \(\chi_k^2\) denotes a chi-square random variable with \(k\) degrees of freedom.

\section{The \texttt{libtube} Library}\label{sec:libtube}

The main computational problem in implementing results based on the
tube formula is evaluation of the constants \(\kappa_0\), \(\kappa_2\),
\(l_0\) e.t.c.
The \texttt{libtube} library implements the tube library up to fourth
order terms. To use the library, one must first write a
`manifold function' defining the problem. \texttt{libtube} takes
the manifold function as input, and uses numerical integration methods
to compute the constants.

The library can be downloaded from
\url{http://www.herine.net/stat/libtube}.
The library is written in C, and can be compiled on Linux systems using
\begin{verbatim}
% make
% make install
\end{verbatim}
to install the libraries in \texttt{/usr/local/lib}.

\begin{verbatim}
% cc -o nlreg nlreg.c -ltube -lmut -lm
\end{verbatim}

The library (and the examples given in this paper) have been written
and tested using the Gnu C compiler available in most Linux distributions.
The C code should be compatible with most other compilers and operating
systems.

\subsection{Manifold Functions}

Suppose the manifold is defined by a vector function
\(T(x)\) mapping a \(d\)-dimensional domain \(\X\) to the
manifold \(\M\) in \(n\)-dimensional space. The constants in
the tube formula can be computed from \(T(x)\) and its derivatives, so
in it's simplest form, the manifold function simply computes these.
In statistical applications, one usually doesn't get \(T(x)\)
naturally, but rather one gets a vector \(l(x)\) such that
\(T(x) = l(x)/\|l(x)\|\) (see the regression examples in
Sections \ref{sec:nlreg} and \ref{sec:scb}). The manifold function
can instead provide \(l(x)\) and its derivatives.

In still other examples, one doesn't even obtain \(l(x)\) directly,
but instead obtains a covariance function
\(\sigma(x,x') = \ip{ l(x)}{l(x')}\) (see the mixture example,
Section \ref{sec:mixture}). Since the distance between any
two points on the manifold is given by
\[
  \|l(x)-l(x')\|^2 = \sigma(x,x) + \sigma(x',x') - 2 \sigma(x,x'),
\]
knowledge of the covariance function determines \(l(x)\) up to an orthogonal
transformation. The manifold function can provide \(\sigma(x,x')\) and
its derivatives.

The precise form of the manifold functions is illustrated by the examples.
After writing the manifold function, the most useful
functions in \texttt{libtube} are:

\begin{itemize}
\item \verb+tube_contstants()+, to numerically evaluate \(\kappa_0\)
  and the other constants appearing in (\ref{eq:}).
\item \verb+tailp()+ and \verb+critval()+,
  which compute tail probabilities corresponding to a specified cut-off,
  and critical values corresponding to a specified significance level.
\end{itemize}

\vspace{3mm}

\noindent\emph{Calling sequence for} \verb+tube_constants()+.

\vspace{3mm}

The function to compute the constants is
\begin{verbatim}
int tube_constants(f, d, n, ev, mg, fl, kap, wk, deb, uc);
int (*f)();
int d, n, ev, mg;
double *fl, *kap, *wk;
int deb, uc;
\end{verbatim}

The arguments to this function are:
\begin{itemize}
\item
\texttt{f} The manifold function to compute \(l(x)\) and its derivatives.
\item
\texttt{d} The dimension of the manifold.
\item
\texttt{m} The maximum length of the \(l(x)\) vectors. The
argument provided is only used to allocate work space; the actual
length of \(l(x)\) is returned by the manifold function.
\item
\texttt{ev} Integration type. For rectangular domains,
  \texttt{ISIMPSON} is the most useful.
\item
\texttt{mg} Integer vector, giving the number of partitions to use
  in each dimension of the numerical integration rules.
\item
\texttt{fl} Integration limits. A numeric vector with length \(2d\).
  The first \(d\) components give lower limits for each variable;
  the remaining \(d\) components give upper limits.
\item
\texttt{kap} is the vector through which the computed constants are returned.
It should be allocated with at least \(\min(d+1,4)\) terms. The values
returned are \(\ko,l_0/2,(\kt+l_1+m_0)/(2\pi)\) and
\((l_2+m_1+n_0)/(4\pi)\).
\item
\texttt{wk} is a workspace vector. If \texttt{wk=NULL}, the required
workspace will be allocated and freed within the \verb+tube_constants()+
function. To pre-allocate the space, the required
length can be found by calling \verb+k0_reqd(d,m)+.
\item
\texttt{terms} Number of terms to compute, from 1 to 4.
\item
\texttt{uc} An indicator variable indicating whether the manifold
  function computes the weight vectors \texttt{uc=0} or
   covariance derivative matrix \texttt{uc=1}.
\end{itemize}

\vspace{3mm}

\noindent\emph{Calling sequence for} \verb+tailp()+
  \emph{and} \verb+critval()+.

\vspace{3mm}

Tail probabilities are computed using the function
\texttt{tailp}:
\begin{verbatim}
double tailp(c,k0,m,d,s,n,process)
double c, *k0, n;
int m, d, s, process;
\end{verbatim}

Critical values corresponding to a specified tail probability are computed
using the \texttt{critval} function:
\begin{verbatim}
double critval(alpha,k0,m,d,s,n,process)
double *k0, al, n;
int m, d, it, s;
\end{verbatim}

These arguments represent:
\begin{itemize}
\item \texttt{c} Cut-off value for \verb+tailp()+.
\item \texttt{alpha} tail probability for \verb+critval()+.
\item \texttt{k0} is the vector of constants computed by
  the \verb+tube_constants()+ function.
\item \texttt{m} is the number of terms in \texttt{k0}. This is the
  returned value of \verb+tube_constants()+, and is equal to
  \(\min(d+1,4)\).
\item \texttt{d} is the dimension of the manifold.
\item \texttt{c} critical value (\texttt{tailp()} only).
\item \texttt{s} Either \verb+ONE_SIDED+ or \verb+TWO_SIDED+.
\item \texttt{n} For the t-process, the residual degrees of freedom
  used to estimate \(\sigma\). For the uniform process, the dimension
  \(n\). Ignored for the Gaussian process. Beware that \(n\) must have
  type double.
\item \texttt{process} Either \texttt{GAUSS} (when \(\epsilon\) is
  multivariate Gaussian); \texttt{TPROC} (Gaussian process with estimated
  variance); or \texttt{UNIF} (when \(\epsilon\) is uniform on the unit
  sphere).
\end{itemize}

\subsection{Writing a Manifold Function with vectors}

The manifold function computes the vector \(l(x)\) and its derivatives.
The basic form of the function is
\begin{verbatim}
int mymf(x,l,reqd)
double *x, *l;
int reqd;
{ /* function body goes here */
}
\end{verbatim}

The \texttt{x} argument is a point in the input space; \texttt{l} is
a vector to be filled in by the manifold function. The final argument,
\texttt{reqd}, is an integer indicating what the library requires from
the manifold function. If \texttt{reqd=0}, only the vector \(l(x)\) is
required. If \texttt{reqd=1}, then both \(l(x)\) and \(l'(x)\) (or all
the first-order partial derivative vectors of \(l(x)\)) are required.
If \texttt{reqd=2}, then additionally the second-order partial derivative
vectors are required.

In statistical applications, the manifold function will generally require
a data vector, sample size \(n\), dimension \(d\)
and variables other than \(x\) in order
to perform its calculations. These variables should be assigned to global
variables so that they are accessible in the manifold function.

The results of the computations are returned through the \texttt{l}
vector. The vector \(l(x)\) is placed in the first \(n\) elements.
The first-order derivatives are placed in the next \(n \times d\)
elements. The second-order derivatives are placed in the next
\(n \times d \times d\) elements.

The function should return \(n\), the length of the vector \(l(x)\).
Generally, this should be equal to the \(n\) value provided in the
\verb+tube_constants()+ call; it should never be larger. It
can be less. An example where it may be less is for a kernel regression
with compactly supported kernel; only the non-zero elements of
\(l(x)\) need be retained.

\subsection{Writing a Manifold Function with a covariance function.}

The structure of a manifold function based on the covariance is
identical to the vector case; it differs in what is computed.

Given a covariance function  \(\sigma(x,x')\), the manifold
function needs to compute (in the one-dimensional case),
\[
  \pmatrix{ \sigma(x,x') &
     \frac{ \partial \sigma(x,x') }{ \partial x'} &
     \frac{ \partial^2 \sigma(x,x') }{ \partial x'^2 } \cr
     \frac{ \partial \sigma(x,x') }{ \partial x} &
     \frac{ \partial^2 \sigma(x,x') }{ \partial x \partial x' } &
     \frac{ \partial^3 \sigma(x,x') }{ \partial x \partial x'^2 } \cr
     \frac{ \partial^2 \sigma(x,x') }{ \partial x^2 } &
     \frac{ \partial^3 \sigma(x,x') }{ \partial x^2 \partial x' } &
     \frac{ \partial^4 \sigma(x,x') }{ \partial x^2 \partial x'^2 } },
\]
evaluated at \(x'=x\). Again, the matrix is stored in the vector
\texttt{l}, with the columns stacked atop each other.

In higher dimensions, the required matrix is most easily written
in terms of differential operators.
The required
\((1+d+d^2) \times (1+d+d^2)\) matrix is
\[
  \pmatrix{ I \cr
    D_{x_1} \cr
     \vdots \cr
    D_{x_d} \cr
    D_{x_1,x_1} \cr
     \vdots \cr
    D_{x_d,x_d} }
    \sigma(x,x')
  \pmatrix{ I &
    D_{x'_1} &
     \ldots &
    D_{x'_d} &
    D_{x'_1,x'_1} &
     \ldots &
    D_{x'_d,x'_d} }
\]
where \(D\) represents the partial derivative operator with respect to
the subscripted variables.

Another view is as follows. If \(\matr{L}\) is the matrix computed by
a manifold function with vectors, then \( \matr{L}^T \matr{L}\) is the
matrix computed by a manifold function with a covariance function.

\section{Example: Testing in Nonlinear Regression}\label{sec:nlreg}

This was the motivating example for \cite{hotelling}, and
was developed in much more detail by \cite{knowsig}.
Suppose one has data \((x_i,Y_i), i =1,\ldots,n\), and a nonlinear
regression model, such as
\begin{equation}
  Y_i = \alpha e^{\gamma x_i} + \epsilon_i. \label{eq:lrmod}
\end{equation}
The important feature of this model is that the parameter \(\alpha\) enters the
model linearly, while \(\gamma\) enters nonlinearly.
Assume that the errors are independent \(N(0,\sigma^2)\).

Consider the problem of testing \(H_0: \alpha=0\) vs \(H_1 : \alpha \ne 0\).
It can be shown that the log-likelihood ratio test statistic is equivalent to
\begin{equation}
  L = \frac{ \min_{\alpha,\gamma} \| Y - a l(\gamma) \|^2 }{ \|Y\|^2 }
  \label{eq:nlrlr}
\end{equation}
where \(l(\gamma)^T = ( e^{\gamma x_1}, \ldots, e^{\gamma x_n})\).

In classical statistical theory, log-likeliood ratio statistics often
have asymptotic \(\chi^2\) distributions. However, this is not the case
for the statistic (\ref{eq:nlrlr}). One way to see this is
to recall that proofs of the \(\chi^2\) results are based on a quadratic
expansion of the statistic under the null parameters. For the present
problem this would require an expansion around \((0,\gamma_0)\) where
\(\gamma_0\) is `the' null value of \(\gamma\). Unfortunately this
is undefined: when \(\alpha=0\), the parameter \(\gamma\) does not
appear in (\ref{eq:lrmod}); it is not identifiable!

For fixed \(\gamma\), minimizing over \(a\) is a linear least-sqaures
problem. It follows that
\[
  L = 1 - \sup_{\gamma} \ip{ \frac{l(\gamma)}{\|l(\gamma)\|} }{
    \frac{Y}{\|Y\|} }^2.
\]
The null hypothesis \(H_0\) is rejected if \(L \le 1-w^2\) for some \(w>0\),
or equivalently, if
\[
  \sup_{\gamma} \left | \ip{ \frac{l(\gamma)}{\|l(\gamma)\|} }{
    \frac{Y}{\|Y\|} } \right | \ge w.
\]
The constant \(w\) must be chosen to obtain a specified significance level.
That is, we need to be able to evaluate probabilities of the form
\begin{equation}
  P( \sup_{\gamma} | \ip{T(\gamma)}{U} | \ge w )
  \label{eq:tgu}
\end{equation}
where \(T(\gamma) = l(\gamma)/ \|l(\gamma)\|\) defines a curve on
the unit sphere, and \(U = Y/\|Y\|\) is (under \(H_0:\alpha=0\))
uniformly distributed on the surface of the sphere.

\subsection{Non-linear Regression: Implementation}

Code implementing the tube formla for the non-linear regression problem
is shown below. The program consists of the manifold function
\texttt{regmf()}, and the main routine \texttt{main()} that reads in
the data and computes the tube constants. Note that the data vectors
and sample size are stored as global variables, so that they can be
accessed within the manifold function.

The manifold function computes the components of \(l(\gamma)\);
\( l_i (\gamma) = e^{\gamma x_i}\), and of \(l'(\gamma)\),
\( l_i'(\gamma) = \gamma e^{\gamma x_i}\). These vectors are
stored end-to-end in the \texttt{l} argument.

\begin{verbatim}
#include <stdio.h>
#include <math.h>
#include <tube.h>
#define MAXN 1000

double x[MAXN], y[MAXN];
int n;

int regmf(gam,l,reqd)
double *gam, *l;
int reqd;
{ int i;
  double *l1;
  l1 = &l[n];
  for (i=0; i<n; i++)
  { l[i] = exp(gam[0]*x[i]);
    l1[i] = x[i]*exp(gam[0]*x[i]);
  }
  return(n);
}

int main()
{ FILE *infile;
  char filename[100];
  int i, mg;
  double gamlimits[2], kappa[4];
  printf("Data filename ? "); scanf("%s",filename);
  printf("n = ? "); scanf("%d",&n);
  infile = fopen(filename,"r");
  for (i=0; i<n; i++) fscanf(infile,"%lf%lf",&x[i],&y[i]);
  gamlimits[0] = -2.0;
  gamlimits[1] = 2.0;
  mg = 100;

  tube_constants(regmf,1,n,ISIMPSON,&mg,gamlimits,kappa,NULL,0,0);
  printf("%8.5f %8.5f\n",kappa[0],kappa[1]);
}
\end{verbatim}

\section{Example: Simultaneous Confidence Bands} \label{sec:scb}

Application of the tube formula to find simultaneous confidence bands for regression
models has been studied in \cite{naiman87}, \cite{sunloascb} among others.
Consider again regression data, but now suppose that the model is
\[
  Y_i = a_0 + a_1 x_i + a_2 x_i^2 + \epsilon_i = \mu(x_i) + \epsilon_i
\]
(although we formulate the problem for quadratic regression, extension
to other linear models is straightforward).
The goal is to find confidence bands
\[
  \hat\mu(x)  \pm c \sqrt{ \hbox{var} ( \hat\mu(x))}
\]
with simultaneous coverage over some nice domain \(\X\):
\begin{equation}
  P( | \hat\mu(x) - \mu(x)| \le c \sigma \|l(x)\| \hbox{ for all } x \in \X)
    = 1-\alpha. \label{eq:scbprob}
\end{equation}

The least-squares estimates of the parameters are
\[
  \pmatrix{ \hat a_0 \cr \hat a_1 \cr \hat a_2 }
    = (\matr{X}^T \matr{X})^{-1} \matr{X}^T Y
\]
where \(\matr{X}\) is the design matrix.
For fixed \(x\), \(\mu(x)\) is estimated by
\[
  \hat\mu(x) = \hat a_0 + \hat a_1 x + \hat a_2 x^2
    = \pmatrix{ 1 & x & x^2 } ( \matr{X}^T \matr{X})^{-1} \matr{X}^T Y
    = \ip{ l(x) }{Y},
\]
and the variance of the estimate is \( \hbox{var}(\hat\mu(x)) = \sigma^2 \|l(x)\| \).

Now, \(\hat\mu(x)-\mu(x) = \ip{l(x)}{\epsilon}\), and the probability
(\ref{eq:scbprob}) is equivalent to
\[
  \alpha = P \paren{ \sup_x | \ip{ \frac{l(x)}{\|l(x)\|} }{ \epsilon }
     | > c  }.
\]
This problem can be solved using the Gaussian process variant of the tube problem.

Suppose \(\matr{X} = \matr{Q}\matr{R}\) is the \(QR\)-decomposition of the
design matrix. Then \(l(x)\) lies in the column space of \(\matr{Q}\) for
all \(x\), and so
\[
  Z(\gamma) = \ip{ \frac{\matr{Q}^T l(x)}{\|\matr{Q}^T l(x)\|} }{
    \matr{Q}^T \epsilon },
\]
so it suffices to work with the vector
\( l^*(x) = \matr{Q}^T l(x) = (\matr{R}^T)^{-1} f(x)\)
where \(f(x)\) is a vector of the polynomial basis functions.
The derivatives are easily found;
\[
  \frac{d}{dx} l^*(x) = (\matr{R}^T)^{-1} \frac{d}{dx} f(x)
\]
and so on.

\subsection{Simultaneous Confidence Bands: Implementation}

Code for the quadratic regression computations, in an arbitrary
number of dimensions, is shown below.
The functions \texttt{quad()}, \texttt{quadi()} and \texttt{quadij()}
compute the quadratic basis functions \(f(x)\); first-order partial
derivatives and second-order partial derivatives respectively.
The manifold function is \texttt{quadmf()}. The \texttt{main()}
function reads the data from a file; computes the design matrix
and its QR-decomposition; and then calls the
\verb+tube_constants()+ function (The QR functions,
\texttt{qr()} and \texttt{qrtinvx()}, as well as \texttt{transpose()},
are part of the \texttt{mut} library).

When the program is run, the user is prompted for a data file (containing
a matrix of the predictor variables); data dimension (\(n\) and \(d\)),
and limits for the confidence band computation.

The tube constants are computed, then the critical value \(c\) for 95\%
confidence bands. Note that the final argument to \texttt{critval} is
the residual degrees of freedom used to estimate \(\sigma\);
\cite{sunloascb} give the modification of (\ref{eq:}) for this case.

\begin{verbatim}
#include <stdio.h>
#include <tube.h>
#include <mutil.h>

int dim, n, p;
double *X;

void quad(x,f)
double *x, *f;
{ int i, j, k;
  k = 0;
  f[k++] = 1.0;
  for (i=0; i<dim; i++) f[k++] = x[i];
  for (i=0; i<dim; i++)
    for (j=i; j<dim; j++)
      f[k++] = x[i]*x[j];
}

void quadi(x,f,i0)
double *x, *f;
int i0;
{ int i, j, k;
  k = 0;
  f[k++] = 0.0;
  for (i=0; i<dim; i++) f[k++] = (i==i0);
  for (i=0; i<dim; i++)
    for (j=i; j<dim; j++)
      f[k++] = (i==i0)*x[j] + (j==i0)*x[i];
}

void quadij(x,f,i0,j0)
double *x, *f;
int i0, j0;
{ int i, j, k;
  k = 0;
  f[k++] = 0.0;
  for (i=0; i<dim; i++) f[k++] = 0.0;
  for (i=0; i<dim; i++)
    for (j=i; j<dim; j++)
      f[k++] = ((i==i0) & (j==j0)) + ((i==j0) & (j==i0));
}

int quadmf(x,l,reqd)
double *x, *l;
int reqd;
{ int i, j, k;
  k = 0;
  quad(x,l);
  qrtinvx(X,l,n,p);
  k++;
  for (i=0; i<dim; i++)
  { quadi(x,&l[k*p],i);
    qrtinvx(X,&l[k*p],n,p);
    k++;
  }
  for (i=0; i<dim; i++)
    for (j=0; j<dim; j++)
    { quadij(x,&l[k*p],i,j);
      qrtinvx(X,&l[k*p],n,p);
      k++;
    }
  return(p);
}

int main()
{ FILE *infile;
  char filename[100];
  int i, j, mg[100];
  double xlim[100], kappa[4], datarow[100];
  printf("Data filename ? "); scanf("%s",filename);
  printf("n = ? "); scanf("%d",&n);
  printf("dim = ? "); scanf("%d",&dim);
  infile = fopen(filename,"r");
  if (infile==NULL)
  { printf("Error: can't read input file\n");
    return(0);
  }
  p = 1 + dim + dim*(dim+1)/2;
  X = (double *)calloc(n*p,sizeof(double));
  for (i=0; i<n; i++)
  { for (j=0; j<dim; j++)
      fscanf(infile,"%lf",&datarow[j]);
    quad(datarow,&X[i*p]);
  }
  transpose(X,n,p);
  qr(X,n,p,NULL);
  for (i=0; i<dim; i++) mg[i] = 20;
  xlim[0] = -2; xlim[1] = -2;
  xlim[2] = 2; xlim[3] = 2;
  tube_constants(quadmf,dim,p,ISIMPSON,mg,xlim,kappa,NULL,3,0);
  printf("kappa: %8.5f %8.5f %8.5f %8.5f\n",kappa[0],kappa[1],kappa[2],kappa[3]);
}
\end{verbatim}

\section{Example: Mixture Models} \label{sec:mixture}

Suppose \(X_1,\ldots,X_n\) are an i.i.d. sample from a density
\[
  f_{\alpha,\lambda}(x) = (1-\alpha)f_0(x) + \alpha \phi(x,\lambda)
\]
where \(\alpha\) and \(\lambda\) are unknown parameters, with
\( 0 \le \alpha \le 1\). The object is to test \(H_0: \alpha=0\)
vs \(H_1: \alpha > 0\). This is a simple example of mixture testing:
under \(H_0\), the single component \(f_0(x)\) describes the data,
while under \(H_1\), the two components are required.

Consider the normalized score process proposed by
\cite{pillloa}. The score process is
\[
  S(\lambda) = \sum_{i=1}^n \frac{ \phi(X_i,\lambda)}{f_0(X_i)} - 1.
\]
Under the null hypothesis, this has mean 0 and covariance function
\(n \sigma(\lambda,\lambda^{\dag})\), where
\[
  \sigma(\lambda,\lambda^{\dag})
    = \int \frac{ \phi(x,\lambda) \phi(x,\lambda^{\dag}) }{f_0(x)} dx - 1.
\]
The normalized score process is
\(S^*(\lambda) = S(\lambda)/\sqrt{n\sigma(\lambda,\lambda)}\).
This asymptotically behaves like a Gaussian
process \(Z(\lambda)\), with mean 0 under \(H_0\), and a nonzero
mean under \(H_1\). The maximum of the normalized score process
serves as the test statistic.

Since an explicit vector representation of \(Z(\lambda)\) is not
readily available, the manifold function (\texttt{mixmf} in the code
below) must be written using the covariance function and its partial 
derivatives.

There is one additional difficulty. The normalized score process
has a singularity at \(\mu=0\). For this reason, the manifold function
works with Taylor series expansions of the covariance in this reason.
Also, the singularity results in a discontinuity in \(S^*(\lambda)\), and
\(l_0=4\).

The main routine in the program below sets limits for \(\mu\), calls the
\verb+tube_constants()+ function, and computes the 5\% critical value.
\begin{verbatim}
kappa0 =  5.27449
  l0/2 =  2.00000
Level 0.05 critical value =  2.49455
\end{verbatim}

\begin{verbatim}
#include <stdio.h>
#include <math.h>
#include <tube.h>
#define MAXN 1000

double x[MAXN], y[MAXN];
int n;

int mixmf(mu,l,reqd)
double *mu, *l;
int reqd;
{ double emm, mm;

  if (fabs(mu[0]) < 0.01)
  { mm = mu[0]*mu[0];
    l[0] = 1 + mm/2*(1 + mm/3*(1 + mm/4*(1 + mm/5)));
    l[1] = 0.5*(1 + mm/3*(2 + mm/4*(3 + mm/5*(4 + 5/6*mm))));
    l[1] = l[2] = mu[0]*l[1];
    l[3] = 0.5*(1 + mm/3*(4 + mm/4*(9 + mm/5*(16+25/6*mm))));
  } else
  { emm = exp(mu[0]*mu[0]);
    l[0] = emm-1;
    l[1] = l[2] = mu[0]*emm;
    l[3] = emm*(1+mu[0]*mu[0]);
  }
  return(2);
}

int main()
{ int i, mg, t;
  double mulimits[2], kappa[4];
  mulimits[0] = -3.0;
  mulimits[1] = 3.0;
  mg = 200;

  t = tube_constants(mixmf,1,n,ISIMPSON,&mg,mulimits,kappa,NULL,2,1);
  /* modify kappa[1] = l0/2 for the singularity */
  kappa[1] += 1.0;

  printf("kappa0 = %8.5f\n",kappa[0]);
  printf("  l0/2 = %8.5f\n",kappa[1]);
  printf("Level 0.05 critical value = %8.5f\n",critval(kappa,t,1,0.05,10,1,0.0));
}
\end{verbatim}

\bibliographystyle{apalike}
\bibliography{paper}

\end{document}